\documentclass[a4, 12pt]{amsart}
\usepackage[mathscr]{eucal}
\usepackage{amssymb}
\usepackage{latexsym}
\usepackage{amsthm}
\theoremstyle{plain}
\newtheorem{theorem}{Theorem}[section]

\newtheorem{remark}{Remark}[section]
\newtheorem{proposition}{Proposition}[section]

\makeatletter
\@addtoreset{equation}{section}

\title[Complete hypersurfaces with $w$-constant mean curvature]{Complete hypersurfaces
with $w$-constant mean curvature  in the unit spheres}
\author{Qing-Ming Cheng and  Guoxin Wei}
\address{Qing-Ming Cheng \\ Department of Applied Mathematics, Faculty of Science,
Fukuoka  University, 814-0180, Fukuoka,  Japan, cheng@fukuoka-u.ac.jp}
\address{Guoxin Wei \\  School of Mathematical Sciences, South China Normal University,
510631, Guangzhou,  China, weiguoxin@tsinghua.org.cn}

\begin{document}
\maketitle

\begin{abstract} In this paper, we study  $4$-dimensional complete hypersurfaces  with $w$-constant mean curvature in the unit sphere. We give a lower  bound
of the scalar curvature  for $4$-dimensional complete hypersurfaces  with $w$-constant mean curvature.
As a  by-product, we give a new  proof of  the result of Deng-Gu-Wei \cite{dgw}  under the  weaker topological condition.
\end{abstract}

\footnotetext{The first author was partially  supported by JSPS Grant-in-Aid for Scientific Research:  No. 22K03303, 
the fund of Fukuoka University: No. 225001.
The second author was partly supported by grant No. 12171164 of NSFC, GDUPS (2018).}

\maketitle
% ----------------------------------------------------------------
\section{ Introduction}
\noindent
It is a very  important problem in differential geometry to investigate the  topology
and geometry of  hypersurfaces with constant mean curvature in spheres.
In order to deal with this  problem,
geometric analysis  is a useful tool.  Especially, for compact minimal hypersurfaces,
by computing  the Laplacian
of the squared norm $S$ of the second fundamental form of
minimal hypersurfaces in spheres, Simons in \cite{s}
proved that  for an $n$-dimensional compact
minimal hypersurface in a unit sphere $S^{n+1}(1)$, if $S\leq n$,
then $S\equiv 0$ or $S \equiv n$.  In the landmark papers  of
Chern, do Carmo and Kobayashi \cite{cdk}  and
Lawson \cite{l}, they proved that  Clifford torus $S^m(\sqrt{\frac mn})\times
S^{n-m}(\sqrt{\frac {n-m}n})$  for $1\leq m\leq n-1$  are the only compact minimal
hypersurfaces in $S^{n+1}(1)$ with  $S\equiv n$.
\\
Furthermore, Chern in \cite{cdk} proposed
the following, which is proposed again by Yau in \cite{y1}:

\vskip3mm
\noindent
{\bf Chern  problems}.
For $n$-dimensional compact minimal  hypersurfaces in $S^{n+1}(1)$
with constant scalar curvature, does the following hold?
\begin{enumerate}
\item $S\leq c(n)$, where $c(n)$ is a universal constant only depending on  $n$.
\item the value $S$ of the squared norm
of the second fundamental forms should be discrete.
\end{enumerate}

\vskip2mm
\noindent
Simons \cite{s}, Chern, do Carmo and Kobayashi \cite{cdk}  and Lawson \cite{l}
obtained  that the first and the second value are 0 and $n$, respectively.
In particular, when $n=3$, Peng and Terng \cite{pt1}
proved that $S\geq 6$ if $S>3$ and $S=6$
if and only if compact minimal  hypersurfaces in $S^{4}(1)$ are isoparametric.
\\
Furthermore, the following Chern conjecture is well-known:

\vskip2mm
\noindent
{\bf Chern  conjecture}.
For $n$-dimensional compact minimal  hypersurfaces in $S^{n+1}(1)$
with constant scalar curvature, if $S>n$, then  $S\geq 2n$.

\noindent
For  the  Chern conjecture, Yang and the first author (\cite{yc1},\cite{yc2},\cite{yc3}) made an important  breakthrough. They proved
if $S>n$, then $S\geq \dfrac{4n}3$. (cf. \cite{c}, \cite{cwt}, \cite{dx}, \cite{gt}, \cite{gxxz}, \cite{pt}, \cite{pt1},  \cite{sy}, \cite{xx}).
For $n=3$, Chern problems are solved affirmatively  (cf. \cite{c1}).

\vskip2mm
\noindent
Minimal hypersurfaces in $S^{n+1}(1)$ are critical points of  the  area functional
and Willmore hypersurfaces  in $S^{n+1}(1)$ are critical points of  the Willmore  functional
$$
\int_M(S-nH^2)^{\frac{n}2}dM.
$$
A hypersurface is called a Willmore minimal hypersurface if it is a critical point of both the area functional
and the Willmore functional.
In \cite{dgw}, Deng, Gu and Wei  studied  the Chern problems for $n=4$
under the additional condition that  $M^4$  is  a $4$-dimensional compact Willmore minimal
hypersurface in $S^{5}(1)$. They proved
\vskip1mm
\noindent
{\bf Theorem DGW}. For a   $4$-dimensional  compact Willmore minimal hypersurface  in $S^5(1)$, if  its scalar curvature is constant,
then, the scalar curvature is non-negative.

\vskip2mm
\noindent
It is also well-known that hypersurfaces  in $S^{n+1}(1)$ with constant mean curvature are critical points of the area
functional for variations which preserving the volume. Compact hypersurfaces  in $S^{n+1}(1)$ with constant mean curvature
are studied extendedly. For examples, Barbosa-do Carmo \cite{bc} proved that a compact hypersurface in $S^{n+1}(1)$
with constant mean curvature is stable if and only if it is isometric to a sphere.
Cheng-Nakagawa \cite{cn} and Xu \cite{x} proved
that for  an $n$-dimensional compact  hypersurface with constant mean curvature $H$ in the unit sphere $S^{n+1}(1)$,
 if the squared length $S$ of the second fundamental form satisfies
 $$
 S \leq n+\dfrac{n^3H^2}{2(n-1)} -\dfrac{n(n-2)}{2(n-1)}\sqrt{ n^2H^4 +4(n-1)H^2},
 $$
then M is either a totally umbilical sphere, or a Clifford torus with constant mean curvature,
which generalized the results  of Chern-do Carmo and Kobayashi \cite{cdk},  Lawson \cite{l} and Simons \cite{s}.
In 1990, de Almeida and Brito \cite{db} proved that compact  3-dimensional hypersurfaces with constant mean curvature
and nonnegative constant scalar curvature in  $S^{n+1}(1)$  must be isoparametric.
Cheng and Wan \cite{cw} and Chang \cite{c} proved , independently,  that a compact hypersurface with constant mean curvature
and constant scalar curvature in the unit sphere $S^{4}(1)$ is  isoparametric. Thus, $3$-dimensional compact hypersurfaces
with constant mean curvature and constant scalar curvature in $S^4(1)$ are classified  completely.
\\
Recently, Tang-Yan\cite{ty}  generalized the theorem of de Almeida and Brito  \cite{db} to higher dimensional cases.
They have proved that for  an  $n$-dimensional compact hypersurface with constant mean curvature and
constant non-negative scalar curvature in  $S^{n+1}(1)$, if  $f_k$, for $k=1, 2, \cdots, n-1$  are constant,
then $M^n$ is isoparametric, where $f_k=\sum_i\lambda_i^k$, $\lambda_i$'s are principal curvatures of the
hypersurfaces.
\vskip2mm
\noindent
On the other hand, Guo and Li \cite{gl} studied the following functional
$$
\int_M(S-nH^2)dM.
$$
For $n=2$,  this functional is a conformal invariant, which is called Willmore functional.  Critical points of  Willmore functional
are  called Willmore surfaces.  Willmore surfaces were  studied by Thomsen, Willmore, Bryant, Pinkall, Weiner, Montiel, Li and
so on.
Guo and Li  \cite{gl} proved that an $n$-dimensional hypersurface $M^n$ in the unit sphere $S^{n+1}(1)$ is a critical point
of the  functional
$$
\int_M(S-nH^2)dM
$$
 if and only if
$$
(n-1)\Delta H=HS+\dfrac{n}2(S-nH^2)H-f_3.
$$
A hypersurface in $S^{n+1}(1)$ is called {\it a hypersurface with $w$-constant mean curvature} in $S^{n+1}(1)$ if it is a critical
point of both the area functional for variations which preserve volume and the functional
$$
\int_M(S-nH^2)dM.
$$
Thus, we know that a hypersurface in $S^{n+1}(1)$ with constant mean curvature is a  hypersurface with
$w$-constant mean curvature  in $S^{n+1}(1)$
if
\begin{equation}\label{eq:1.1}
f_3=HS+\dfrac{n}2(S-nH^2)H.
\end{equation}
In this paper, we study  complete hypersurfaces with $w$-constant mean curvature  in $S^{n+1}(1)$. We will prove the following:
\begin{theorem}\label{theorem 1.1}
For a $4$-dimensional  complete hypersurface  $M^4$  with $w$-constant mean curvature in $S^5(1)$, the scalar curvature $R$
satisfies 
 $$
 R> -3.3(1+H^2)
 $$
  if  the  scalar curvature  $R$ is constant.
\end{theorem}
\noindent
As a by-product, we prove the following:

\begin{theorem}\label{theorem 1.2}
For a   $4$-dimensional  complete Willmore minimal  hypersurface  $M^4$  in $S^5(1)$,  the scalar curvature is non-negative
 if  the scalar curvature is constant.
\end{theorem}
\begin{remark}
In our theorem, the topological condition is weaker than  one in Deng-Gu-Wei \cite{dgw} and our proof is  different from one in \cite{dgw}.
\end{remark}
\vskip2mm
\noindent

\begin{remark}
We propose that  a   $4$-dimensional  complete  hypersurface  $M^4$ in $S^5(1)$ with $w$-constant mean curvature  and  constant
scalar curvature has non-negative scalar curvature, that is,
$$
S\leq 12(1+H^2)+4H^2.
$$
In order to prove this assertion, we think  that  all informations on the second covariant derivative $h_{ijkl}$ of
the  second fundamental form are necessary.  But its computation is  very complicated.
\end{remark}

\section{Preliminary}
\vskip2mm
\noindent
Let $M^n$ be an $n$-dimensional
hypersurface in $S^{n+1}(1)$. We choose a local orthonormal frame
$\{\vec{e}_1, \cdots, \vec{e}_{n}, \vec{e}_{n+1}\}$ and the dual
coframe $\{\omega_1, \cdots,$ $\omega_n$, $ \omega_{n+1}\}$ in
such  that $\{\vec{e}_1, \cdots, \vec{e}_n\}$  is a local
orthonormal frame on $M^n$. Hence, we have
$$
\omega_{n+1}=0
$$
in $M^n$. According to Cartan lemma, one has
$$
\omega_{i,n+1}=\sum_{j}h_{ij}\omega_j, \ h_{ij}=h_{ji}.
$$
The mean curvature $H$ and the second fundamental form
$\vec{\alpha}$ of $M^n$ are defined, respectively, by
$$
H=\frac{1}{n}\sum_{i}h_{ii}, \
\vec{\alpha}=\sum_{i,j}h_{ij}\omega_i\otimes\omega_j\vec{e}_{n+1}.
$$
When  $H\equiv 0$,
 $M^n$ is   called a  minimal hypersurface.
From the structure equations of $M^n$,  we have Gauss equations,
and Codazzi equations.
$$
\begin{aligned}
&R_{ijkl}=(\delta_{ik}\delta_{jl}-\delta_{il}\delta_{jk})
+(h_{ik}h_{jl}-h_{il}h_{jk}),\\
\end{aligned}
$$
\begin{align*}
h_{ijk}=h_{ikj},
\end{align*}
where $h_{ijk}$ are defined by
\begin{equation}
\sum_{k}h_{ijk}\omega_k=dh_{ij}+\sum_k
h_{ik}\omega_{kj} +\sum_k h_{kj}\omega_{ki},
\end{equation}
Defining $h_{ijkl}$, $h_{ijklm}$  by
\begin{equation}
\sum_{l}h_{ijkl}\omega_l=dh_{ijk}+\sum_l
h_{ijl}\omega_{lk} +\sum_l h_{ilk}\omega_{lj}+\sum_lh_{ljk}\omega_{li},
\end{equation}
\begin{equation}
\begin{aligned}
&\sum_{m}h_{ijklm}\omega_m\\&=dh_{ijkl}+\sum_m
h_{ijkm}\omega_{ml} +\sum_m h_{ijml}\omega_{mk}\sum_m h_{imkl}\omega_{mj} +\sum_m h_{mjkl}\omega_{mi},
\end{aligned}
\end{equation}
we have  Ricci formulas
\begin{align*}
h_{ijkl}-h_{ijlk}=\sum_mh_{im}R_{mjkl}+\sum_mh_{mj}R_{mikl}.
\end{align*}
\begin{align*}
h_{ijklm}-h_{ijkml}=\sum_ph_{ijp}R_{pklm}+\sum_ph_{ipk}R_{pjlm}+\sum_ph_{pjk}R_{pilm}.
\end{align*}
We define functions $f_3$ and $f_4$ by
\begin{align*}
f_3=\sum_{i,j,k}h_{ij}h_{jk}h_{ki} \ \ {\rm and} \ \
f_4=\sum_{i,j,k,l}h_{ij}h_{jk}h_{kl}h_{li}
\end{align*}
respectively.
At each point, we may take a local frame $\{\vec{e}_1, \cdots, \vec{e}_n\}$ such that
$$
h_{ij}=\lambda_i\delta_{ij}.
$$
One has
$$
\begin{aligned}
&h_{ij}-H\delta_{ij}=\lambda_i\delta_{ij}-H\delta_{ij}=:\mu_i\delta_{ij}, \ B:=\sum_{i}\mu_{i}^2=S-nH^2, \\
& B_3=\sum_{i}\mu_i^3, \  B_4=\sum_{i}\mu_i^4, \ \  B_5=\sum_{i}\mu_i^5,
\end{aligned}
$$
$$
R_{ijkl}=\bigl(1+H^2+H(\mu_i+\mu_j)+\mu_i\mu_j\bigl)(\delta_{ik}\delta_{jl}-\delta_{il}\delta_{jk}).
$$
By a simple computation, we have  from Ricci formulas
\begin{proposition} For a hypersurface with constant mean curvature in $S^{n+1}(1)$, we have
\begin{equation}
\begin{aligned}
\Delta h_{ij}&=(n(1+H^2)-B)(h_{ij}-H\delta_{ij})+H(n\mu_i^2-B)\delta_{ij},\\
\Delta h_{ijk}&=(2n+3+2nH^2-B)h_{ijk}+nH(\mu_i+\mu_j+\mu_k)h_{ijk}\\
&-(\lambda_i^2+\lambda_j^2+\lambda_k^2)h_{ijk}-(\nabla_kBh_{ij}+\nabla_jBh_{ik}+\nabla_iBh_{jk})\\
&+2(\lambda_k\lambda_i+\lambda_j\lambda_k+\lambda_i\lambda_j)h_{ijk}.
\end{aligned}
\end{equation}
\end{proposition}

\noindent
Thus, we obtain
\begin{proposition} For a hypersurface with constant mean curvature
in $S^{n+1}(1)$, we have
\begin{equation}\label{eq:2.5}
\begin{aligned}
&\dfrac12\Delta B=\sum_{i,j,k}h_{ijk}^2+(n(1+H^2)-B)B+nHB_3,\\
&\dfrac12\Delta \sum_{i,j,k}h_{ijk}^2=\sum_{i,j,k,l}h_{ijkl}^2+(2n+3+2nH^2-B)\sum_{i,j,k}h_{ijk}^2+3nH\sum_{i,j,k}\mu_ih_{ijk}^2\\
&\qquad\ \ \ \ \ \ \ \ \ \ \ -3\sum_{i,j,k}\lambda_i^2h_{ijk}^2+6\sum_{i,j,k}\lambda_k\lambda_ih_{ijk}^2-\dfrac32|\nabla B|^2,\\
&\dfrac13\Delta B_3=(n+nH^2-B)B_3+nHB_4-HB^2+2\sum_{i,j,k}\mu_ih_{ijk}^2,\\
& \dfrac14\Delta B_4= (n+nH^2-B)B_4+H(nB_5-BB_3) +2\sum_{i,j,k}\mu_i^2h_{ijk}^2+\sum_{i,j,k}\mu_i\mu_jh_{ijk}^2.
\end{aligned}
\end{equation}
\end{proposition}
\noindent
From Gauss equations, we know
\begin{equation}
R=n(n-1)+n^2H^2-S=n(n-1)(1+H^2)-B.
\end{equation}
Hence, for  a hypersurface with constant mean curvature
in $S^{n+1}(1)$, scalar curvature $R$ is constant if and only if $B$ is constant.
\\
The following Generalized Maximum Principle due to Omori \cite{o} (cf. Yau \cite{y}) will play an important role in this paper.
\begin{theorem}
Let $M^{n}$ be a complete Riemannian manifold with sectional  curvature bounded from below.
If a $\mathcal C^2$-function $f$ is bounded from above in $M^{n}$, then there exists a sequence of $\{p_{k}\}_{k=1}^{\infty} \subset{M^{n}}$ such that
\begin{enumerate}
\item $\lim_{k \to \infty} f(p_{k})=\underset{M^{n}}{\sup}\ f$,
\item $\lim_{k \to \infty}|\nabla f(p_{k})|=0$,
\item $\lim_{k \to \infty}\sup\ \nabla_l\nabla_lf(p_{k})\leq 0$, for  $l=1,  2,  \cdots, n$.
\end{enumerate}
\end{theorem}

\section{ Proofs of theorems}

\vskip2mm
\noindent
For a hypersurface with $w$-constant mean curvature, we have
$$
f_3=HS+\dfrac{n}2(S-nH^2)H.
$$
If $n=4$, we have
\begin{equation}\label{eq:3.0}
B_3=\sum_i\mu_i^3=f_3-3HB-4H^3=0, \quad \sum_{i,j,k}h_{ijk}^2=(B-4(1+H^2))B.
\end{equation}
For convenient, we may assume
$$
\mu_1\geq \mu_2\geq \mu_3\geq\mu_4.
$$
According to
$$
B_3=(\mu_1+\mu_4)(\mu_1^2+\mu_4^2-\mu_1\mu_4-\mu_2^2-\mu_3^2+\mu_2\mu_3),
$$
we have
$$
\mu_1+\mu_4=0, \quad  \mu_2+\mu_3=0.
$$
\vskip2mm
\noindent
{\it Proof of theorem 1.1}.
Since the scalar curvature is constant, we know that $S$ is constant. According to Gauss equations, we
know that the sectional curvature of $M^4$ is bounded. By applying the Generalized Maximum Principle  due to
Omori \cite{o} (cf. Yau \cite{y}) to the function $-B_4$,
 there exists a sequence of $\{p_{k}\}_{k=1}^{\infty}\subset M^{4}$ such that
\begin{equation*}
\lim_{k \to \infty}B_{4}(p_{k})=\inf\ B_{4},\
\lim_{k \to \infty}|\nabla B_{4}(p_{k})|=0,\
\lim_{k \to \infty}\sup\nabla_l\nabla_l B_{4}(p_{k})\geq 0, \text{\rm for} \ l=1, 2, 3, 4.
\end{equation*}
Since $H$, $S$ are  constant,  we have from the proposition 2.2
\begin{equation}\label{eq:3.1}
\begin{aligned}
&\sum_{i,j,k}h_{ijk}^2=(B-n(1+H^2))B,\\
&\sum_{i,j,k,l}h_{ijkl}^2=(B-2n-3-2nH^2)\sum_{i,j,k}h_{ijk}^2-3nH\sum_{i,j,k}\mu_ih_{ijk}^2\\
&\ \ \ \ \ \ \ \ \qquad+3\sum_{i,j,k}\lambda_i^2h_{ijk}^2-6\sum_{i,j,k}\lambda_k\lambda_ih_{ijk}^2.\\
\end{aligned}
\end{equation}
Hence, we know that, for any $i ,j ,k ,l$,   $\{\lambda_{i}(p_{k})\}$,  \ $\{h_{ijk}(p_{k})\}$ and $\{h_{ijkl}(p_{k})\}$ are bounded sequences, respectively.
Thus, we can assume, if necessary,  by taking  a subsequences of $\{p_{m}\}$,
\begin{equation*}
\lim_{m \to \infty} \lambda_{i}(p_{{m}})=\hat{\lambda}_{i},\
\lim_{m \to \infty} h_{ijk}(p_{{m}})=\hat{h}_{ijk},\
\lim_{m \to \infty} h_{ijkl}(p_{{m}})=\hat{h}_{ijkl},\ \ \forall i,j,k,l.
\end{equation*}
From now on, for  all of  computations on $\hat{\lambda}_{i}, \hat{h}_{ijk}$ and $\hat{h}_{ijkl}$ ,  we omit \ $\hat{}$ \  for convenience.\newline
Since $\mu_1+\mu_4=-(\mu_2+\mu_3)=0$,  if $ \mu_1=\mu_2$, by making use the almost same assertions as in \cite{dgw}, we can  prove that 
$M^4$ is isometric to the  Clifford torus with constant mean curvature. We would like to omit the detailed proof in this case.\\
If three of  $\mu_1, \mu_2, \mu_3, \mu_4$ are equal to each other, we know $\mu_1=-\mu_4>0$ and $ \mu_2= \mu_3=0$,
$$
 \inf B_4=\frac {B^2}2=\sup B_4.
 $$
We know that $B_4=\dfrac {B^2}4$ is constant and $M^4$ has three distinct constant principal curvatures  at each point.
It is impossible.

 \vskip2mm
 \noindent
From now, we can assume 
$$
\mu_1>\mu_2>\mu_3>\mu_4,  \ \  \mu_1+\mu_4=-\mu_2-\mu_3=0.
$$
Since $H$, $B$ , $B_3=0$ are  constant, we obtain under our limiting process
\begin{equation}
\begin{cases}
&\sum_ih_{iik}=0, \\
&\sum_i\mu_ih_{iik}=0,\\
&\sum_i\mu_i^2h_{iik}=0,\\
&\sum_i\mu_i^3h_{iik}=0.\\
\end{cases}
\end{equation}
Hence,  we have
$$
h_{iik}=0,
$$
for any $i, k=1, 2, 3, 4$. \newline
From  $B_3=0$,
we have
$$
nHB_4-HB^2+2\sum_{i,j,k}\mu_ih_{ijk}^2=0.
$$
$\lim_{k \to \infty}|\nabla B_{4}(p_{k})|=0$ and
$\sum_{i,j,k}h_{ijk}^2=B(B-4-4H^2)$,
we conclude, for any $m$,
$$
\begin{aligned}
&\sum_{i,j,k}h_{ijk}h_{ijkm}=0,\\
&2\sum_{i,j,k}h_{ijkm}h_{ijk}\mu_i+\sum_{i,j,k}h_{ijk}h_{ijl}h_{klm}=0.
\end{aligned}
$$
$$
\begin{aligned}
&h_{123}h_{123m}+h_{124}h_{124m}+h_{134}h_{134m}+h_{234}h_{234m}=0,\\
&\mu_1h_{123}h_{123m}+\mu_2h_{124}h_{124m}+\mu_3h_{134}h_{134m}+\mu_4h_{234}h_{234m}\\
&=-\bigl(h_{123}h_{134}h_{24m}+h_{123}h_{124}h_{34m}+h_{123}h_{234}h_{14m}\\
&\ \ \ +h_{124}h_{234}h_{13m}+h_{124}h_{134}h_{23m}+h_{134}h_{234}h_{12m}\bigl).
\end{aligned}
$$
Because $B_3=0$ is constant, by a direct computation, we have
\begin{equation}\label{eq:3.2}
\sum_{i,j,k}\mu_i^2h_{ijk}^2-2\sum_{i,j,k}\mu_i\mu_jh_{ijk}^2=BB_4-(1+H^2)B^2.
\end{equation}
By taking covariant derivative of the above equality and taking limiting process, we obtain
$$
\sum_{i,j,k}(\mu_k^2-2\mu_i\mu_j)h_{ijk}h_{ijkm}=\sum_{i,j,k}(\mu_i+\mu_j-\mu_k)h_{ijk}h_{ijl}h_{klm},
$$
$$
\begin{aligned}
&\mu_1^2h_{123}h_{123m}+\mu_2^2h_{124}h_{124m}+\mu_3^2h_{134}h_{134m}+\mu_4^2h_{234}h_{234m}\\
&=(\mu_2+\mu_4)h_{123}h_{134}h_{24m}+(\mu_3+\mu_4)h_{123}h_{124}h_{34m}\\
&\ \ \ +(\mu_1+\mu_3)h_{124}h_{234}h_{13m}+(\mu_1+\mu_2)h_{134}h_{234}h_{12m}.
\end{aligned}
$$
Hence, we obtain
\begin{equation}
\begin{cases}
&h_{123}h_{123m}+h_{124}h_{124m}+h_{134}h_{134m}+h_{234}h_{234m}=0,\\
&\mu_1h_{123}h_{123m}+\mu_2h_{124}h_{124m}+\mu_3h_{134}h_{134m}+\mu_4h_{234}h_{234m}\\
&=-\bigl(h_{123}h_{134}h_{24m}+h_{123}h_{124}h_{34m}+h_{123}h_{234}h_{14m}\\
&\ \ +h_{124}h_{234}h_{13m}+h_{124}h_{134}h_{23m}+h_{134}h_{234}h_{12m}\bigl),\\
&\mu_1^2h_{123}h_{123m}+\mu_2^2h_{124}h_{124m}+\mu_3^2h_{134}h_{134m}+\mu_4^2h_{234}h_{234m}\\
&=(\mu_2+\mu_4)h_{123}h_{134}h_{24m}+(\mu_3+\mu_4)h_{123}h_{124}h_{34m}\\
&\ \ +(\mu_1+\mu_3)h_{124}h_{234}h_{13m}+(\mu_1+\mu_2)h_{134}h_{234}h_{12m}.
\end{cases}
\end{equation}
By solving the above linear equation systems, we have
\begin{equation}\label{eq:3.4}
\begin{cases}
\mu_2(\mu_2^2-\mu_1^2)h_{123}h_{1231}&=\mu_2(\mu_1^2-\mu_2^2)h_{234}h_{1234}+2\mu_1\mu_2h_{123}h_{124}h_{134},\\
\mu_2(\mu_2^2-\mu_1^2)h_{124}h_{1241}&=-\mu_1(\mu_1^2-\mu_2^2)h_{234}h_{1234}\\
&\ \ \ -\dfrac12(\mu_1+\mu_2)(3\mu_2-\mu_1)h_{123}h_{124}h_{134},\\
\mu_2(\mu_2^2-\mu_1^2)h_{134}h_{1341}&=-\mu_1(\mu_1^2-\mu_2^2)h_{234}h_{1234}\\
&\ \ \ -\dfrac12(\mu_1-\mu_2)(3\mu_2+\mu_1)h_{123}h_{124}h_{134},\\
\end{cases}
\end{equation}
\begin{equation}\label{eq:3.5}
\begin{cases}
\mu_1(\mu_1^2-\mu_2^2)h_{123}h_{1232}&=\mu_2(\mu_1^2-\mu_2^2)h_{134}h_{1234}\\
&\ \ \ -\frac12(\mu_1+\mu_2)(3\mu_1-\mu_2)h_{123}h_{124}h_{234},\\
\mu_1(\mu_1^2-\mu_2^2)h_{124}h_{1242}&=-\mu_1(\mu_1^2-\mu_2^2)h_{134}h_{1234}
+2\mu_1\mu_2h_{123}h_{124}h_{234},\\
\mu_1(\mu_1^2-\mu_2^2)h_{134}h_{2342}&=-\mu_2(\mu_1^2-\mu_2^2)h_{134}h_{1234}\\
&\ \ \ +\frac12(\mu_1-\mu_2)(3\mu_2+\mu_1)h_{123}h_{124}h_{234},\\
\end{cases}
\end{equation}
\begin{equation}\label{eq:3.6}
\begin{cases}
\mu_1(\mu_1^2-\mu_2^2)h_{123}h_{1233}&=-\mu_2(\mu_1^2-\mu_2^2)h_{124}h_{1234}\\
&\ \ \ -\frac12(\mu_1-\mu_2)(3\mu_1+\mu_2)h_{123}h_{134}h_{234},\\
\mu_1(\mu_1^2-\mu_2^2)h_{134}h_{1343}&=-\mu_1(\mu_1^2-\mu_2^2)h_{124}h_{1234}
-2\mu_1\mu_2h_{123}h_{134}h_{234},\\
\mu_1(\mu_1^2-\mu_2^2)h_{234}h_{2343}&=\mu_2(\mu_1^2-\mu_2^2)h_{124}h_{1234}\\
&\ \ \ +\frac12(\mu_1+\mu_2)(3\mu_1-\mu_2)h_{123}h_{134}h_{234},\\
\end{cases}
\end{equation}
\begin{equation}\label{eq:3.7}
\begin{cases}
\mu_2(\mu_1^2-\mu_2^2)h_{124}h_{1244}&=-\mu_1(\mu_1^2-\mu_2^2)h_{123}h_{1234}\\
&\ \ \ -\frac12(\mu_1-\mu_2)(\mu_1+3\mu_2)h_{124}h_{134}h_{234},\\
\mu_2(\mu_1^2-\mu_2^2)h_{134}h_{1344}&=\mu_1(\mu_1^2-\mu_2^2)h_{123}h_{1234}\\
&\ \ \ +\frac12(\mu_1+\mu_2)(\mu_1-3\mu_2)h_{124}h_{134}h_{234},\\
\mu_2(\mu_1^2-\mu_2^2)h_{234}h_{2344}&=-\mu_2(\mu_1^2-\mu_2^2)h_{123}h_{1234}
+2\mu_1\mu_2h_{124}h_{134}h_{234}.\\
\end{cases}
\end{equation}
Because  $B_3=0$ and $H$ is constant , by taking covariant derivative, we have, for any $k, m$,
$$
\sum_{i}h_{iikm}=0, \  \ \sum_{i}\mu_i^2h_{iikm}+\sum_{i,j}(\mu_i+\mu_j)h_{ijk}h_{ijm}=0.
$$
Since $\sum_{i,j}(\mu_i+\mu_j)h_{ij2}h_{ij3}=0$ and $\sum_{i,j}(\mu_i+\mu_j)h_{ij1}h_{ij4}=0$, we obtain
\begin{equation}
\begin{aligned}
(\mu_1^2-\mu_2^2)(h_{1123}+h_{4423})=0, \  \ (\mu_1^2-\mu_2^2)(h_{2214}+h_{3314})=0,
\end{aligned}
\end{equation}
that is, $h_{1123}+h_{4423}=0$ and $h_{2214}+h_{3314}=0$.\newline
We get
$(h_{123}^2+h_{234}^2)h_{1234}=0$ and $(h_{124}^2+h_{134}^2)h_{1234}=0$.
Since $h_{123}^2+h_{234}^2 +h_{124}^2+h_{134}^2\neq 0$, we have $h_{1234}=0$.\newline
Since $\sum_{i,j}(\mu_i+\mu_j)h_{ij1}h_{ij2}=-2(\mu_1+\mu_2)h_{134}h_{234}$ and
$$
\sum_{i}h_{ii12}=0, \  \ \sum_{i}\mu_i^2h_{ii12}+\sum_{i,j}(\mu_i+\mu_j)h_{ij1}h_{ij2}=0,
$$
we obtain
$$
h_{1112}+h_{4412}=2\dfrac{h_{134}h_{234}}{\mu_1-\mu_2},
\ \ \ h_{2212}+h_{3312}=-2\dfrac{h_{134}h_{234}}{\mu_1-\mu_2}.
$$
Because $B$ is constant, we know
$$
\sum_{i}\mu_ih_{iikm}+\sum_{i,j}h_{ijk}h_{ijm}=0, \ \ \ \ \ \ \forall k,m.
$$
From $\sum_{i,j}h_{ij1}h_{ij2}=2h_{134}h_{234}$, we have
$$
\mu_1(h_{1112}-h_{4412})+\mu_2(h_{2212}-h_{3312})=-2h_{134}h_{234}.
$$
Hence, we infer
\begin{equation}\label{eq:6-28-1}
\mu_1h_{4412}+\mu_2h_{3312}=2h_{134}h_{234}.
\end{equation}
According to $h_{1234}=0$,  \eqref{eq:3.6}  and \eqref{eq:3.7}, we have
$$
\begin{aligned}
2\mu_1(\mu_1+\mu_2)h_{123}h_{1233}&=-(3\mu_1+\mu_2)h_{123}h_{134}h_{234},\\
2\mu_2(\mu_1+\mu_2)h_{124}h_{1244}&=-(\mu_1+3\mu_2)h_{124}h_{134}h_{234}.\\
\end{aligned}
$$
Hence
\begin{equation}\label{eq:6-28-2}
h_{123}h_{124}(\mu_1h_{1233}+\mu_2h_{1244})=-2h_{123}h_{124}h_{134}h_{234}.
\end{equation}
From \eqref{eq:6-28-1} and \eqref{eq:6-28-2}, we have
$$
\begin{aligned}
&h_{123}h_{124}(\mu_1-\mu_2)h_{1244}=2h_{123}h_{124}h_{134}h_{234},\\
&h_{123}h_{124}(\mu_1-\mu_2)h_{1233}=-2h_{123}h_{124}h_{134}h_{234},\\
\end{aligned}
$$
$$
\begin{aligned}
4\mu_2(\mu_1+\mu_2)h_{123}h_{124}h_{134}h_{234}&=-(\mu_1-\mu_2)(\mu_1+3\mu_2)h_{123}h_{124}h_{134}h_{234},\\
\end{aligned}
$$
namely,
$$
\bigl\{(\mu_1+\mu_2)^2+4\mu_1\mu_2\bigl\}h_{123}h_{124}h_{134}h_{234}=0.
$$
We obtain, from $\mu_1>\mu_2>0$,
$$
h_{123}h_{124}h_{134}h_{234}=0.
$$
We can assume $h_{123}=0$ or $h_{134}=0$.\newline
Since $H$, $B$ , $B_3=0$ are  constant, we obtain under our limiting process, for any $k, l$,
\begin{equation}
\begin{cases}
&\sum_ih_{iikl}=0, \\
&\sum_i\mu_ih_{iikl}+\sum_{i,j}h_{ijk}h_{ijl}=0,\\
&\sum_i\mu_i^2h_{iikl}+\sum_{i,j}(\mu_i+\mu_j)h_{ijk}h_{ijl}=0.\\
\end{cases}
\end{equation}
Hence, we have
\begin{equation}
\begin{cases}
&2\mu_2h_{22kk}=-2\mu_1h_{11kk}-\dfrac{\sum_{i,j}(\mu_i+\mu_j)h_{ijk}^2}{\mu_1+\mu_2}-\sum_{i,j}h_{ijk}^2, \\
&2\mu_2h_{33kk}=2\mu_1h_{11kk}+\mu_1\dfrac{\sum_{i,j}(\mu_i+\mu_j)h_{ijk}^2}{\mu_1^2-\mu_2^2}+\sum_{i,j}h_{ijk}^2, \\
&h_{44kk}=-h_{11kk}-\dfrac{\sum_{i,j}(\mu_i+\mu_j)h_{ijk}^2}{\mu_1^2-\mu_2^2},\\
\end{cases}
\end{equation}
\begin{equation}
\begin{aligned}
&\dfrac14\nabla_k\nabla_kB_4=\sum_{i}\mu_i^3h_{iikk}+2\sum_{i,j}\mu_i^2h_{ijk}^2+\sum_{i,j}\mu_i\mu_jh_{ijk}^2\\
&=\mu_1(\mu_1^2-\mu_2^2)(h_{11kk}-h_{44kk})-\mu_2^2\sum_{i,j}h_{ijk}^2+2\sum_{i,j}\mu_i^2h_{ijk}^2+\sum_{i,j}\mu_i\mu_jh_{ijk}^2\\
&=2\mu_1(\mu_1^2-\mu_2^2)h_{11kk}+\mu_1\sum_{i,j}(\mu_i+\mu_j)h_{ijk}^2\\
&\ \ \ -\mu_2^2\sum_{i,j}h_{ijk}^2+2\sum_{i,j}\mu_i^2h_{ijk}^2+\sum_{i,j}\mu_i\mu_jh_{ijk}^2.\\
\end{aligned}
\end{equation}

\noindent
According to \eqref{eq:2.5},  \eqref{eq:3.0}, \eqref{eq:3.2} and $h_{iik}=0$ for any $i, k$, we have
\begin{equation}
\begin{cases}
&h_{123}^2+h_{124}^2+h_{134}^2+h_{234}^2=\dfrac16B(B-4-4H^2),\\
&\mu_4h_{123}^2+\mu_3h_{124}^2+\mu_2h_{134}^2+\mu_1h_{234}^2= H(B_4-\dfrac14B^2),\\
&\mu_4^2h_{123}^2+\mu_3^2h_{124}^2+\mu_2^2h_{134}^2+\mu_1^2h_{234}^2=\dfrac1{9}B^2(B-4-4H^2)\\
&\qquad \ \ \ \ \ \ \ \ \ \ \ \ \ \ \ \ \ \ \ \ \ \ \ \ \ \ \ \ \ \ \ \ \ \ \ \ \ \ \ \ \ -\dfrac16\bigl\{BB_4-(1+H^2)B^2\bigl\}.\\
\end{cases}
\end{equation}
\begin{equation}
\begin{aligned}
\dfrac14\Delta B_4&=\dfrac83(B-6-6H^2)\mu_1^2\mu_2^2-\dfrac{B^2}{18}(B-10-10H^2).\\
\end{aligned}
\end{equation}
If $h_{134}=0$, we have, from $B_4=\dfrac{B^2}4+(\mu_1^2-\mu_2^2)^2$,

\begin{equation}
\begin{cases}
&h_{123}^2+h_{124}^2+h_{234}^2=\dfrac16B(B-4-4H^2),\\
&\mu_4h_{123}^2-\mu_2h_{124}^2+\mu_1h_{234}^2=H(\mu_1^2-\mu_2^2)^2,\\
&\mu_4^2h_{123}^2+\mu_2^2h_{124}^2+\mu_1^2h_{234}^2
=\dfrac5{72}B^2(B-4-4H^2)-\dfrac B6(\mu_1^2-\mu_2^2)^2.\\
\end{cases}
\end{equation}
By  solving this linear equation system, we have
\begin{equation}
\begin{aligned}
2\mu_1(\mu_1-\mu_2)h_{123}^2&=\dfrac{5B^2}{72}(B-4-4H^2)-\dfrac{B}6(\mu_1^2-\mu_2^2)^2\\
&\ \ \ -\dfrac{B}{6}(B-4-4H^2)\mu_1\mu_2-H(\mu_1-\mu_2)(\mu_1^2-\mu_2^2)^2,\\
(\mu_1^2-\mu_2^2)h_{124}^2&=\dfrac{B}6\bigl[(\mu_1^2-\dfrac{5B}{12})(B-4-4H^2)+(\mu_1^2-\mu_2^2)^2\bigl],\\
2\mu_1(\mu_1+\mu_2)h_{234}^2&=\dfrac{5B^2}{72}(B-4-4H^2)-\dfrac{B}6(\mu_1^2-\mu_2^2)^2\\
&\ \ \ +\dfrac{B}{6}(B-4-4H^2)\mu_1\mu_2+H(\mu_1+\mu_2)(\mu_1^2-\mu_2^2)^2.\\
\end{aligned}
\end{equation}
From (3.10),  we have, by using Ricci formulas,
\begin{equation}
\begin{cases}
h_{2211}=&-\dfrac{\mu_1}{\mu_2}h_{1111}-\dfrac{1}{\mu_2}h_{123}^2-\dfrac{2}{\mu_1+\mu_2}h_{124}^2,\\
h_{3311}=&\dfrac{\mu_1}{\mu_2}h_{1111}+\dfrac{1}{\mu_2}h_{123}^2,\\
h_{4411}=&- h_{1111}+\dfrac{2}{\mu_1+\mu_2}h_{124}^2,\\
\end{cases}
\end{equation}

\begin{equation}
\begin{cases}
h_{1144}=&-h_{1111}+\dfrac{2}{\mu_1+\mu_2}h_{124}^2+2\mu_1(1+H^2-\mu_1^2),\\
h_{2244}=&\dfrac{\mu_1}{\mu_2}h_{1111}-\dfrac{2(2\mu_1+\mu_2)}{\mu_2(\mu_1+\mu_2)}h_{124}^2
-\dfrac1{\mu_2}h_{234}^2-2\dfrac{\mu_1^2}{\mu_2}(1+H^2-\mu_1^2),\\
h_{3344}=&-\dfrac{\mu_1}{\mu_2}h_{1111}+\dfrac{4\mu_1^2}{\mu_2(\mu_1^2-\mu_2^2)}h_{124}^2
+\dfrac1{\mu_2}h_{234}^2+2\dfrac{\mu_1^2}{\mu_2}(1+H^2-\mu_1^2),\\
h_{4444}=&h_{1111}-\dfrac{4\mu_1}{\mu_1^2-\mu_2^2}h_{124}^2-2\mu_1(1+H^2-\mu_1^2).\\
\end{cases}
\end{equation}

Hence, we have
\begin{equation}\label{eq:3.20}
\begin{aligned}
&\dfrac14\nabla_1\nabla_1B_4=2\mu_1(\mu_1^2-\mu_2^2)h_{1111},\\
&\dfrac14\nabla_4\nabla_4B_4=-2\mu_1(\mu_1^2-\mu_2^2)h_{1111}+8\mu_1^2h_{124}^2
+4\mu_1^2(\mu_1^2-\mu_2^2)(1+H^2-\mu_1^2),\\
&\dfrac14\nabla_1\nabla_4B_4=-Bh_{123}h_{234},\\
&\dfrac14\Delta B_4=\dfrac83(B-6-6H^2)\mu_1^2\mu_2^2-\dfrac{B^2}{18}(B-10-10H^2).
\end{aligned}
\end{equation}
Putting $\mu_1^2=t\mu_2^2$ and $a=\dfrac{1+H^2}{\mu_2^2}$, we have
$$
t>1, \ \ a<\dfrac{1+t}6,
$$
if $B>12(1+H^2)$.
$$
\begin{aligned}
&2\mu_1(\mu_1-\mu_2)h_{123}^2\\
&=\biggl\{\dfrac{2(1+t)}{9}\bigl[t^2-3t\sqrt t+8t-3\sqrt t +1-(5t-6\sqrt t+5)a\bigl]
-\dfrac{H}{\mu_2}(\sqrt t-1)(t-1)^2\biggl\}\mu_2^6,\\
\end{aligned}
$$
$$
\begin{aligned}
&2\mu_1(\mu_1+\mu_2)h_{234}^2\\
&=\biggl\{\dfrac{2(1+t)}{9}\bigl[t^2+3t\sqrt t+8t+3\sqrt t +1-(5t+6\sqrt t+5)a\bigl]
+\dfrac{H}{\mu_2}(\sqrt t+1)(t-1)^2\biggl\}\mu_2^6,\\
\end{aligned}
$$
$$
\begin{aligned}
&(\mu_1^2-\mu_2^2)h_{124}^2=\dfrac{2(1+t)}{9}\biggl\{2t^2-5t-1-(t-5)a\biggl\}\mu_2^6.
\end{aligned}
$$
Since
\begin{equation*}
\begin{aligned}
&\dfrac14\nabla_1\nabla_1B_4=2\mu_1(\mu_1^2-\mu_2^2)h_{1111},\\
&\dfrac14\nabla_4\nabla_4B_4=-2\mu_1(\mu_1^2-\mu_2^2)h_{1111}
+8\mu_1^2h_{124}^2
+4\mu_1^2(\mu_1^2-\mu_2^2)(1+H^2-\mu_1^2),\\
&\dfrac14\nabla_1\nabla_4B_4=-Bh_{123}h_{234},\\
&\dfrac14\Delta B_4=\dfrac83(B-6-6H^2)\mu_1^2\mu_2^2-\dfrac{B^2}{18}(B-10-10H^2).
\end{aligned}
\end{equation*}
Thus, we have
\begin{equation}\label{eq:3.21}
\begin{aligned}
\dfrac14\Delta B_4&=\biggl\{\dfrac83(2(1+t)-6a)t-\dfrac{2(1+t)^2}9(2(1+t)-10a)\biggl\}\mu_2^6\\
&=\biggl\{\dfrac{4(t+1)}{9}\bigl\{12\bigl(1-3\dfrac{a}{t+1}\bigl)t-(1+t)^2(1-5\dfrac{a}{t+1})\bigl\}\biggl\}\mu_2^6\\
&=\biggl\{\dfrac{4(t+1)}{9}\bigl\{-t^2+10t-1+(5t^2-26t+5)\dfrac{a}{t+1}\bigl\}\biggl\}\mu_2^6,\\
\end{aligned}
\end{equation}
\begin{equation}\label{eq:3.22}
\begin{aligned}
&\dfrac14\nabla_1\nabla_1B_4+\dfrac14\nabla_4\nabla_4B_4=8\mu_1^2h_{124}^2
+4\mu_1^2(\mu_1^2-\mu_2^2)(1+H^2-\mu_1^2)\\
&=\biggl\{\dfrac{16t(1+t)}{9(t-1)}\bigl\{2t^2-5t-1-(t-5)a\bigl\}+4t(t-1)(a-t)\biggl\}\mu_2^6\\
&=\biggl\{\dfrac{4(t+1)}{9}\bigl\{-t^2+6t-34-\dfrac{16}{t-1}+\dfrac{18}{t+1}\\
&\ \ +(5t^2+3t+32+\dfrac{32}{t-1})\dfrac{a}{t+1}\bigl\}\biggl\}\mu_2^6.\\
\end{aligned}
\end{equation}
If $h_{123}=0$,
we have
\begin{equation}
\begin{cases}
&h_{124}^2+h_{134}^2+h_{234}^2=\dfrac16B(B-4-4H^2),\\
&-\mu_2h_{124}^2+\mu_2h_{134}^2+\mu_1h_{234}^2= H(B_4-\dfrac14B^2),\\
&\mu_2^2h_{124}^2+\mu_2^2h_{134}^2+\mu_1^2h_{234}^2=\dfrac1{9}B^2(B-4-4H^2)\\
&\ \ \ \ \ \ \ \ \ \ \ \ \ \ \ \ \ \ \ \ \ \ \ \ \ \ \ \ \ \ \ \ \ \ \ -\dfrac16\bigl\{BB_4-(1+H^2)B^2\bigl\}.\\
\end{cases}
\end{equation}
\begin{equation}
\begin{aligned}
&(\mu_1^2-\mu_2^2)h_{234}^2\\
&=\dfrac19(B-\dfrac{3}{2}\mu_2^2)B(B-4-4H^2)-\dfrac16\bigl\{BB_4-(1+H^2)B^2\bigl\}\\
&=(\dfrac{5}{72}B-\dfrac{1}{6}\mu_2^2)B(B-4-4H^2)-\dfrac{B}{6}(\mu_1^2-\mu_2^2)^2\\
&=\dfrac{2(t+1)}{9}\bigl\{t^2+5t-2-(5t-1)a\bigl\}\mu_2^6,
\end{aligned}
\end{equation}
\begin{equation}\label{eq:3.25}
\begin{aligned}
&2\mu_2h_{134}^2\\
&=-(\mu_1+\mu_2)h_{234}^2+\dfrac{\mu_2}6B(B-4-4H^2)+H(B_4-\dfrac14B^2)\\
&=-\dfrac{2(t+1)}{9(\sqrt t-1)}\bigl\{t^2+5t-2-(5t-1)a\bigl\}\mu_2^5\\
&\ \ \ +\dfrac{2(t+1)}{3}(t+1-2a)\mu_2^5
+\dfrac{H}{\mu_2}(t-1)^2\mu_2^5\\
&=-\dfrac{2(t+1)}{9(\sqrt t-1)}\bigl\{(t^2+8t+1)-3\sqrt t(t+1)-(5t+5-6\sqrt t))a\bigl\}\mu_2^5\\
&\ \ \ +\dfrac{H}{\mu_2}(t-1)^2\mu_2^5,
\end{aligned}
\end{equation}
\begin{equation}
\begin{aligned}
&2\mu_2h_{124}^2=(\mu_1-\mu_2)h_{234}^2+\dfrac{\mu_2}6B(B-4-4H^2)-H(B_4-\dfrac14B^2)\\
&=\dfrac{2(t+1)}{9(\sqrt t+1)}\bigl\{t^2+5t-2-(5t-1)a\bigl\}\mu_2^5\\
&\ \ \ +\dfrac{2(t+1)}{3}(t+1-2a)\mu_2^5-\dfrac{H}{\mu_2}(t-1)^2\mu_2^5\\
&=\dfrac{2(t+1)}{9(\sqrt t+1)}\bigl\{(t^2+8t+1+3\sqrt t(t+1)-(5t+5+6\sqrt t)a\bigl\}\mu_2^5\\
&\ \ \ -\dfrac{H}{\mu_2}(t-1)^2\mu_2^5.
\end{aligned}
\end{equation}
\vskip2mm
\noindent
From (3.10) and (3.11), $h_{iik}=0$, $h_{123}=0$, we obtain
\begin{equation}
\begin{aligned}
&\dfrac14\nabla_1\nabla_1B_4=2\mu_1(\mu_1^2-\mu_2^2)h_{1111},\\
\end{aligned}
\end{equation}
\begin{equation}
\begin{aligned}
&\dfrac14\nabla_4\nabla_4B_4=-2\mu_1(\mu_1^2-\mu_2^2)h_{1111}+8\mu_1^2h_{124}^2+8\mu_1^2h_{134}^2\\
&\ \ \ \ \ \ \ \ \ \ \ \ \ \ \ \ +4\mu_1^2(\mu_1^2-\mu_2^2)(1+H^2-\mu_1^2)\\
&=-2\mu_1(\mu_1^2-\mu_2^2)h_{1111}+4t(t-1)(a-t)\mu_2^6\\
&\ \ \ +\dfrac{8t(t+1)}{9(t-1)}\biggl\{(\sqrt t-1)\bigl\{(t^2+8t+1+3\sqrt t(t+1)-(5t+5+6\sqrt t)a\bigl\}\\
&\ \ \ -(\sqrt t+1)\bigl\{(t^2+8t+1)-3\sqrt t(t+1)-(5t+5-6\sqrt t))a\bigl\}\biggl\}\mu_2^6\\
&=-2\mu_1(\mu_1^2-\mu_2^2)h_{1111}+4t(t-1)(a-t)\mu_2^6\\
&\ \ \ +\dfrac{16t(t+1)}{9(t-1)}\bigl\{2t^2-5t-1-(t-5)a\bigl\}\mu_2^6.\\
\end{aligned}
\end{equation}
Hence, we have
\begin{equation}\label{eq:3.29}
\begin{aligned}
&\dfrac14\nabla_1\nabla_1B_4+\dfrac14\nabla_4\nabla_4B_4\\
&=4t(t-1)(a-t)\mu_2^6
+\dfrac{16t(t+1)}{9(t-1)}\bigl\{2t^2-5t-1-(t-5)a\bigl\}\mu_2^6.\\
\end{aligned}
\end{equation}
Therefore, we have, from \eqref{eq:3.21}, \eqref{eq:3.22} and \eqref{eq:3.29},
if $h_{123}=0$ or $h_{134}=0$,

\begin{equation}\label{eq:3.30}
\begin{aligned}
\dfrac14\Delta B_4&=\dfrac83(B-6-6H^2)\mu_1^2\mu_2^2-\dfrac{B^2}{18}(B-10-10H^2)\\
&=\dfrac{4(t+1)}{9}\bigl\{-t^2+10t-1+(5t^2-26t+5)\dfrac{a}{t+1}\bigl\}\mu_2^6,\\
\end{aligned}
\end{equation}
\begin{equation}\label{eq:3.31}
\begin{aligned}
&\dfrac14\nabla_1\nabla_1B_4+\dfrac14\nabla_4\nabla_4B_4\\
&=4t(t-1)(a-t)\mu_2^6
+\dfrac{16t(t+1)}{9(t-1)}\bigl\{2t^2-5t-1-(t-5)a\bigl\}\mu_2^6.\\
\end{aligned}
\end{equation}
Since
$$
\lim_{k\to\infty}\inf\nabla_l\nabla_lB_4(p_k)\geq  0,
$$
for any $l$, we have,  under our limiting process,

$$
4t(t-1)(a-t)
+\dfrac{16t(t+1)}{9(t-1)}\bigl\{2t^2-5t-1-(t-5)a\geq 0,
$$
$$
\dfrac83(2(1+t)-6a)t-\dfrac{2(1+t)^2}9(2(1+t)-10a)\geq 0,
$$
that is,
\begin{equation}\label{eq:3.33-1}
\begin{aligned}
& -t^2+6t-34-\dfrac{16}{t-1}+\dfrac{18}{t+1}+(5t^2+3t+32+\dfrac{32}{t-1})\dfrac{a}{t+1}\geq 0, \\
&(1+t)(10t-1-t^2)+(5t^2-26t+5)a\geq 0.
\end{aligned}
\end{equation}
From \eqref{eq:3.33-1} and  $\dfrac{a}{t+1}<\dfrac16$, we have
$$
\begin{aligned}
&-t^2+39t-172-\dfrac{64}{t-1}+\dfrac{108}{t+1}\geq 0.\\
\end{aligned}
$$
Defining
$$
f(t)=t^2-39t+172+\dfrac{64}{t-1}-\dfrac{108}{t+1},
$$
we know 
$$
\dfrac{df(t)}{dt}=2t-39-\dfrac{64}{(t-1)^2}+\dfrac{108}{(t+1)^2}<0, \ 1<t<5.
$$
Hence, $f(t)$ is a decreasing function. Since $f(5)=0$, we know that $f(t)>0$ for $t<5$.
This is a contradiction.
Hence, we have $t\geq 5$ and 
$5t^2-26t+5\geq 0$.\\
From  \eqref{eq:3.33-1} and $\dfrac{a}{t+1}<\dfrac16$, we have
$$
-t^2+34t-1\geq 0.
$$
Therefore, we get 
$$
5\leq t\leq 17+12\sqrt 2.
$$
From \eqref{eq:3.33-1}, we have
$$
\begin{aligned}
&-t^2+6t-34-\dfrac{16}{t-1}+\dfrac{18}{t+1}+(5t^2+3t+32+\dfrac{32}{t-1})\dfrac{a}{t+1}\geq 0.\\
\end{aligned}
$$
Hence, we obtain
$$
\begin{aligned}
&\dfrac{a}{t+1}\geq \dfrac{t^2-6t+34+\dfrac{16}{t-1}-\dfrac{18}{t+1}}{(5t^2+3t+32+\dfrac{32}{t-1})}\geq 0.130729.
\\
\end{aligned}
$$
Since 
$B=2(1+t)\dfrac{1+H^2}{a}$, we conclude
$$
B\leq \dfrac{2}{0.130729}(1+H^2)<15.3(1+H^2).
$$
Thus,  from Gauss equation, one has
$$
R>-3.3(1+H^2).
$$
We finish our proof of the theorem \ref{theorem 1.1}.\\

\noindent {\it Proof of Theorem \ref{theorem 1.2}}.

\noindent Since our assertion in proof of the theorem 1.1 is true for $H=0$,
if  $h_{123}=0$, we obtain, from \eqref{eq:3.25},
\begin{equation}
\begin{aligned}
&2\mu_2h_{134}^2\\
&=-\dfrac{2(t+1)}{9(\sqrt t-1)}\bigl\{(t^2+8t+1)-3\sqrt t(t+1)-(5t+5-6\sqrt t))a\bigl\}\mu_2^5.\\
\end{aligned}
\end{equation}
This is impossible since $5\leq t\leq 17+12\sqrt 2$ and $\dfrac{a}{t+1}<\dfrac16$.\\
If $h_{134}=0$,
since
\begin{equation*}
\begin{aligned}
&4\mu_1^2(\mu_1^2-\mu_2^2)h_{123}^2h_{234}^2\\
&=\biggl\{\dfrac{2(1+t)}{9}\bigl[t^2-3t\sqrt t+8t-3\sqrt t +1-(5t-6\sqrt t+5)a\bigl]
-\dfrac{H}{\mu_2}(\sqrt t-1)(t-1)^2\biggl\}\mu_2^6\\
&\ \ \ \times\biggl\{\dfrac{2(1+t)}{9}\bigl[t^2+3t\sqrt t+8t+3\sqrt t +1-(5t+6\sqrt t+5)a\bigl]
+\dfrac{H}{\mu_2}(\sqrt t+1)(t-1)^2\biggl\}\mu_2^6,\\
\end{aligned}
\end{equation*}
that is,
\begin{equation}
\begin{aligned}
&\dfrac{4\mu_1^2(\mu_1^2-\mu_2^2)h_{123}^2h_{234}^2}{\mu_2^{12}}\\
&=\dfrac{4(t+1)^2}{81}\biggl\{\bigl[t^2+8t +1-5(t+1)a\bigl]^2-9t(t+1-2a)^2\biggl\}\\
&\ \ \ -\dfrac{H}{\mu_2}\dfrac{4(1+t)(t-1)^2}{9}\bigl[2t^2-5t-1-( t-5)a\bigl]-\dfrac{H^2}{\mu_2^2}(t-1)^5,
\end{aligned}
\end{equation}
according to

$$
\dfrac1{16}(\Delta B_4)^2\geq \dfrac1{16}\{\nabla_1\nabla_1B_4+\nabla_4\nabla_4B_4\}^2
\geq \dfrac1{4}\nabla_1\nabla_1B_4\cdot\nabla_4\nabla_4B_4\geq \dfrac1{4}(\nabla_1\nabla_4B_4)^2,
$$
we have, from \eqref{eq:3.20} and \eqref{eq:3.21},

\begin{equation}
\begin{aligned}
&\dfrac{16(t+1)^2}{81}\bigl\{-t^2+10t-1+(5t^2-26t+5)\dfrac{a}{t+1}\bigl\}^2\\
&\geq \dfrac{16(t+1)^4}{81t(t-1)}
\biggl\{\bigl[t^2+8t +1-5(t+1)a\bigl]^2-9t(t+1-2a)^2\biggl\},\\
\end{aligned}
\end{equation}
that is,
\begin{equation}
\begin{aligned}
&t(t-1)\bigl\{t^2-10t+1-(5t^2-26t+5)\dfrac{a}{t+1}\bigl\}^2\\
&\geq (t+1)^2
\biggl\{\bigl[t^2+8t +1-5(t+1)a\bigl]^2-9t(t+1-2a)^2\biggl\}\\
&=(t+1)^2\biggl\{\bigl[t^2-10t+1-(5t^2-26t+5)\dfrac{a}{t+1}\bigl]^2\\
&\ \ \ +36t(1-2\dfrac{a}{t+1})(t^2-10t+1-(5t^2-26t+5)\dfrac{a}{t+1})\\
&\ \ \ +(18)^2t^2(1-2\dfrac{a}{t+1})^2
-9t(t+1)^2(1-2\dfrac{a}{t+1})^2\biggl\}\\
&=(t+1)^2\biggl\{\bigl[t^2-10t+1-(5t^2-26t+5)\dfrac{a}{t+1}\bigl]^2\\
&\ \ \ +27t(1-2\dfrac{a}{t+1})(t^2-2t+1)(1-6\dfrac{a}{t+1})\biggl\}.\\
\end{aligned}
\end{equation}
Hence, we obtain
\begin{equation}
\begin{aligned}
&-(3t+1)\bigl\{t^2-10t+1-(5t^2-26t+5)\dfrac{a}{t+1}\bigl\}^2\\
&\geq 27t(t-1)^2(1-2\dfrac{a}{t+1})(1-6\dfrac{a}{t+1})>0.
\end{aligned}
\end{equation}
This is a  contradiction. We must have $S\leq 12$.  From the Gauss equations, we know that
the scalar curvature is non-negative.
This completes the proof of Theorem \ref{theorem 1.2}.\\

% ----------------------------------------------------------------
\bibliographystyle{amsplain}

\end {document}